\documentclass[twoside]{article} 
\usepackage{graphicx}
\title{\bf Asymptotics of the Mittag-Leffler function $E_a(z)$ on the negative real axis when $a\to1$}
\author{\sc R. B. Paris\footnote{E-mail address:\ \ {\tt r.paris@abertay.ac.uk}}\\
\\
{\em Division of Computing and Mathematics,}\\
{\em Abertay University, Dundee DD1 1HG, UK}\\
}

\begin{document}
\newcommand{\bee}{\begin{equation}}
\newcommand{\ee}{\end{equation}}
\def\f#1#2{\mbox{${\textstyle \frac{#1}{#2}}$}}
\def\dfrac#1#2{\displaystyle{\frac{#1}{#2}}}
\newcommand{\fr}{\frac{1}{2}}
\newcommand{\fs}{\f{1}{2}}
\newcommand{\g}{\Gamma}
\newcommand{\al}{\alpha}
\newcommand{\bb}{\beta}
\newcommand{\ta}{\theta}
\newcommand{\om}{\omega}
\newcommand{\br}{\biggr}
\newcommand{\bl}{\biggl}
\newcommand{\ra}{\rightarrow}
\renewcommand{\topfraction}{0.9}
\renewcommand{\bottomfraction}{0.9}
\renewcommand{\textfraction}{0.05}
\newcommand{\mcol}{\multicolumn}
\newcommand{\gtwid}{\raisebox{-.8ex}{\mbox{$\stackrel{\textstyle >}{\sim}$}}}
\newcommand{\ltwid}{\raisebox{-.8ex}{\mbox{$\stackrel{\textstyle <}{\sim}$}}}
\date{}
\maketitle
\pagestyle{myheadings}
\markboth{\hfill {\it R.B. Paris} \hfill}
{\hfill {\it The Mittag-Leffler function $E_a(z)$ as $a\to 1$} \hfill}
\begin{abstract} 
We consider the asymptotic expansion of the single-parameter Mittag-Leffler function $E_a(-x)$ for $x\to+\infty$ as the parameter $a\to1$. The dominant expansion when $0<a<1$ consists of an algebraic expansion of $O(x^{-1})$ (which vanishes when $a=1$), together with an exponentially small contribution that approaches $e^{-x}$ as $a\to 1$. Here we concentrate on the form of this exponentially small expansion when $a$ approaches the value 1.

Numerical examples are presented to illustrate the accuracy of the expansion so obtained. 
\vspace{0.4cm}

\vspace{0.3cm}

\end{abstract}

\vspace{0.2cm}

\noindent $\,$\hrulefill $\,$

\vspace{0.2cm}

\begin{center}
{\bf 1. \  Introduction}
\end{center}
\setcounter{section}{1}
\setcounter{equation}{0}
\renewcommand{\theequation}{\arabic{section}.\arabic{equation}}
The single-parameter Mittag-Leffler function $E_a(z)$ is defined by
\bee\label{e11}
E_a(z)=\sum_{n=0}^\infty \frac{z^n}{\Gamma(an+1)} \qquad (|z|<\infty)
\ee
where $a>0$. This function has recently found application in fractional calculus  and in the modelling of `non-standard' processes; see, for example, \cite{GKMR, HMS, MG, SR}. When $0<a<1$, it also arises 
in the standard model of fractional diffusion \cite{MLP}. In particular, when $z=-x$ ($x>0$), the limit $a\to1$  corresponds to the transition from fractional (slow) diffusion to classical diffusion.

In this paper we shall restrict the parameter $a$ to satisfy $0<a<1$ and pay particular attention to the above-mentioned limit $a\to 1$. 
The standard asymptotic expansion of $E_a(z)$ for $|z|\to\infty$ when $0<a<1$ is \cite[\S 18.1]{Erd}, \cite[\S 5.1.4]{PK}
\bee\label{e12}
E_a(z)\sim\left\{\begin{array}{ll}\dfrac{1}{a} \exp (z^{1/a})+H(z) & (|\arg\,z|<\pi a),\\ 
\\
H(z) & (|\arg (-z)|<\pi(1-\fs a)), \end{array}\right.
\ee  
where the algebraic expansion $H(z)$ is given by the formal asymptotic sum
\bee\label{e13}
H(z)=-\sum_{k=1}^\infty \frac{z^{-k}}{\g(1-ak)}=-\frac{1}{\pi}\sum_{k=1}^\infty \g(ak) \sin (\pi ak)\,z^{-k}.
\ee
When $a=1$, $H(z)\equiv 0$ and the Mittag-Leffler function reduces to the simple exponential function $e^z$.

In the first expansion in (\ref{e12}) we have extended the domain of validity of the compound expansion up to the Stokes lines $\arg\,z=\pm\pi a$. In the sector $|\arg\,z|<\fs\pi a$, the exponential term is dominant for large $|z|$, becoming oscillatory in character on $\arg\,z=\pm\fs\pi a$. In the sectors $\fs\pi a<|\arg\,z|<\pi a$ the exponential term is subdominant and, although exponentially small, can still make a significant contribution in high-precision asymptotics.
On the rays $\arg\,z=\pm\pi a$, the exponential term is maximally subdominant relative to the algebraic expansion $H(z)$. Across these rays a Stokes phenomenon occurs, where in the sense of increasing $|\arg\,z|$ the exponential term ``switches off'' in a smooth manner described approximately by an error function of appropriate argument \cite{MVB},
\cite[p.~67]{DLMF}. A detailed discussion of the Stokes phenomenon associated with $E_a(z)$ for $a>0$ both in $\arg\,z$ (at fixed $|z|$) and in the parameter $a$ is given in \cite{P}. 

Mathematically, we have two distinct asymptotic behaviours, namely $E_a(-x)=O(x^{-1})$ when $0<a<1$ and $E_a(-x)=e^{-x}$ when $a=1$. It is seen from the previous paragraph that the limit $a\to 1$ corresponds to the two Stokes lines $\arg\,z=\pm\pi a$ coalescing on the negative real axis.
The problem of investigating how this transition between dominant algebraic behaviour and the single exponential $e^{-x}$ occurs as $a\to1$ was suggested to the author by W. McLean \cite{McL}. 

It is the purpose of this paper to examine in detail the asymptotic expansion of the Mittag-Leffler function $E_a(-x)$ on the negative real axis as the parameter $a\to 1$. In order for the paper to be self-contained we repeat part of the analysis described in \cite{P} which is given in the appendix. This enables the determination of the form of the exponentially small expansion in this limit and how the single exponential $e^{-x}$ makes its appearance. This is our main result and is stated in Theorem 1 at the beginning of Section 3. 
We supply numerical results in Section 4 that confirm the accuracy of our expansion.
\vspace{0.6cm}

\begin{center}
{\bf 2. \ The expansion of $E_a(-x)$ when $0<a<1$}
\end{center}
\setcounter{section}{2}
\setcounter{equation}{0}
\renewcommand{\theequation}{\arabic{section}.\arabic{equation}}
The two-parameter Mittag-Leffler function $E_{a,b}(z)$ satisfies the recursion property
\[E_{a,b}(z)=\sum_{n=0}^\infty\frac{z^n}{\g(an+b)}=z^{-1} E_{a,b-a}(z)-\frac{z^{-1}}{\g(b-a)}.\]
Application of this result $M$ times, where $M$ is an arbitrary positive integer, yields the result for $E_{a,1}(z)\equiv E_a(z)$
\bee\label{e21}
E_a(z)=-\sum_{k=1}^M\frac{z^{-k}}{\g(1-ak)}+R_M(a;z),\qquad R_M(a;z)=z^{-M} E_{a,1-aM}(z).
\ee
The finite sum on the right-hand side of (\ref{e21}) corresponds to the first $M$ terms of the asymptotic expansion $H(z)$ in (\ref{e13}). We put $z=xe^{i\ta}$, where it is sufficient to consider $0\leq\ta\leq\pi$ since $E_a(xe^{-i\ta})$ is given by the conjugate value.

We shall choose $M$ to be the optimal truncation index of $H(z)$ (corresponding to truncation at, or near, the least term in modulus) given by $aM\sim |z|^{1/a}$ as $|z|\to\infty$. More specifically, we set
\bee\label{e22}
aM=X+\nu,\qquad X=x^{1/a},
\ee
where $\nu$ is bounded. From \cite[\S 18.1]{Erd}, \cite[(2.4)]{P}, we have the integral representation
\bee\label{e22a}
R_M(a;z)=\frac{z^{-M}}{2\pi i} \int_{C'}\frac{u^{a+aM-1} e^u}{u^a-z}\,du=
\frac{e^{-iM\ta}}{2\pi i}\int_C \frac{\tau^{a+aM-1}}{\tau^a-e^{i\ta}}\,e^{X\tau}d\tau,
\ee 
where $C$ denotes a loop surrounding the unit disc with endpoints at $-\infty$ on either side of the branch cut along the negative $\tau$-axis (with $C'$ being the map of this loop in the $u$-plane). The integrand has poles at the points $P_k=\exp\,[i(\ta+2\pi k)/a]$ ($k=0, \pm 1, \pm 2,\ldots$)
and, since $aM\sim X$, it also has saddle points at $e^{\pm\pi i}$; see Fig.~1.
When $\ta<\pi a$, the pole $P_0$ is situated in $0\leq\arg\,\tau<\pi$. The contour $C$ can be deformed over $P_0$ and round the branch point at $\tau=0$ (which is integrable) to yield
the expansion given in \cite[(2.9)]{P}. This expansion, however, breaks down in the vicinity of $\ta=\pi a$ (a Stokes line) since the pole $P_0$ becomes coincident with the saddle point at $\tau=e^{\pi i}$ in this limit.

\begin{figure}[t]
	\begin{center}	\includegraphics[width=0.5\textwidth]{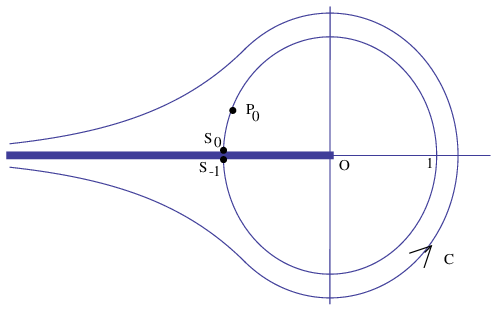}
\caption{\small{The loop $C$ in the $\tau$-plane cut along the negative real axis with the pole $P_0$ and the saddle points $S_0$ and $S_{-1}$ at $e^{\pm\pi i}$. The figure shows the pole corresponding to $\theta<\pi a$; when $\theta=\pi$, $P_0$ crosses the branch cut and passes onto the adjacent Riemann sheet.
}}
	\end{center}
\end{figure}

In what follows we consider the expansion of the remainder term $R_M(a;z)$ in a region enclosing the Stokes line $\ta=\pi a$, which will enable us to deal with the case $\ta=\pi$. For $0\leq\ta\leq\pi a$, the pole $P_0$ lies on the principal Riemann sheet ($|\arg\,\tau|\leq\pi$) in the $\tau$-plane; when $\ta=\pi a$ the pole $P_0$ lies on $\arg\,\tau=\pi$ and when $\ta=\pi$ the pole $P_0$ has passed onto the adjacent sheet that connects with the principal sheet along $\arg\,\tau=\pi$. The pole $P_{-1}=e^{i(\ta-2\pi)/a}$ lies on the adjacent sheet that connects with the principal sheet along $\arg\,\tau=-\pi$. We observe that when $\ta=\pi$, the poles $P_0$ and $P_{-1}$ are situated symmetrically at $e^{\pm\pi i/a}$ on the (separate) adjacent Riemann sheets.   

The loop $C$ in (\ref{e22}) is now deformed round the branch point $\tau=0$, with the path on the upper side of the cut passing {\it above} the pole $P_0$ and saddle at $\tau=e^{\pi i}$ and the path on the lower side of the cut passing {\it below} the pole $P_{-1}$ and saddle at $\tau=e^{-\pi i}$. We concentrate on the contribution from the integral taken along the upper side of the cut since, when $\ta=\pi$, that from the integral along the lower side of the cut will yield the conjugate value.  The details of this calculation are given in the appendix.
\vspace{0.3cm}

\begin{center}
{\bf 3. \ The expansion when $\ta=\pi$}
\end{center}
\setcounter{section}{3}
\setcounter{equation}{0}
\renewcommand{\theequation}{\arabic{section}.\arabic{equation}}
From (\ref{e28}), the expansion of the integral along the upper side of the branch cut in the $\tau$-plane when
$\ta=\pi$ is 
\[\frac{e^{Xe^{\pi i/a}}}{2a}\,\mbox{erfc}\,[c(\pi) \sqrt{X/2}\,]-\frac{ie^{-X-i\om(\pi) X}}{a\sqrt{2\pi X}}\sum_{k=0}^\infty B_{2k}(\pi) (\fs)_k (\fs X)^{-k}\]
as $X\to+\infty$
with $\om(\pi)=\pi(1-a)/a$. The contribution to the integral (\ref{e22a}) from the lower side of the branch cut in the $\tau$-plane will yield the conjugate of the above expansion.

Hence, on the negative real $z$-axis we have the following result:
\newtheorem{theorem}{Theorem}
\begin{theorem}$\!\!\!.$\ \  
The expansion of the Mittag-Leffler function $E_a(-x)$ for $x\to+\infty$ and $\f{1}{3}<a<1$ is
\bee\label{e29}
E_a(-x)= -\sum_{k=1}^M \frac{(-x)^{-k}}{\g(1-ak)}+R_M(a;-x),
\ee
where the remainder $R_M(a;-x)$ has the the exponentially small expansion
\bee\label{e210}
R_M(a;-x)\sim\frac{2}{a} \Re \bl\{\frac{\exp [Xe^{\pi i/a}]}{2}\,\mbox{erfc}\,[c(\pi) \sqrt{ X/2}\,]-\frac{ie^{-X-i\om(\pi)X}}{\sqrt{2\pi X}} \sum_{k=0}^\infty B_{2k}(\pi) (\fs)_k (\fs X)^{-k}\br\}.
\ee  
Here, $M$ is the optimal truncation index of the algebraic expansion given in (\ref{e22}), $X=x^{1/a}$, $\om(\pi)=\pi(1-a)/a$ and $c(\pi)$ is determined from (\ref{e26}) and, in the limit $a\to 1$, by (\ref{e27}). 
\end{theorem}
It now remains to discuss the coefficients $B_{2k}(\pi)$ appearing in (\ref{e210}), which is carried out in the next sub-section.
\vspace{0.5cm}

\noindent 3.1\ \ {\bf The coefficients $B_{2k}(\pi)$}
\vspace{0.3cm}

\noindent
If we use the Series command in {\it Mathematica}, we can obtain the coefficients in the expansion
$f(u)=\sum_{r=0}^\infty \al_ru^r$,
where $f(u)$ is defined in (\ref{e25b}). Upon inversion of the transformation in (\ref{e24}) to obtain
\[t-1=u+\f{1}{3}u^2+\f{1}{36}u^3-\f{1}{270}u^4+\f{1}{4320}u^5+\ldots\,,\]
the first three even-order coefficients $\alpha_{r}$ are found to be
\begin{eqnarray*}
\al_0&=&\frac{1}{1-T},\\
\al_2&=&\frac{1}{12(1-T)^3}\bl\{1+6\nu^2(1-T)^2+(6a^2+6a-2)T+(6a^2-6a+1)T^2\\
&&\hspace{8cm}-6\nu(1-T)(1+(2a-1)T)\br\},\\
\al_4&=&
\frac{1}{864 (1 - T)^5}\bl\{1 + 36 \nu^4 (1 - T)^4 + 
   4 (-1 + 9 a + 30 a^2 + 30 a^3 + 9 a^4) T  \\
   &&
   +6 (1 - 18 a - 20 a^2 + 60 a^3 + 66 a^4) T^2 + 
   4 (-1 + 27 a - 30 a^2 - 90 a^3 + 99 a^4) T^3 \\
   &&+ (1 - 36 a + 
      120 a^2 - 120 a^3 + 36 a^4) T^4 -
   24 \nu^3 (1 - T)^3 (5 + (-5 + 6 a) T) \\
   &&+ 24 \nu^2 (1 - 
      T)^2 (5 + (-10 + 15 a + 9 a^2) T + (5 - 15 a + 9 a^2) T^2)\\
      && - 
   12 \nu (1 - 
      T) (3 + (-9 + 20 a + 30 a^2 + 12 a^3) T + (9 - 40 a + 
         48 a^3) T^2 \\
         &&+ (-3 + 20 a - 30 a^2 + 12 a^3) T^3)\br\},
\end{eqnarray*}
where we have put $T=e^{ia\om(\pi)}$ for brevity.
We recall that the quantity $\nu$ appears in the definition of the optimal truncation index $M$ in (\ref{e22}).
It is impractical to present higher coefficients as they depend on three quantities ($a$, $\nu$ and $T$) and rapidly become too complicated. However, when dealing with specific cases, where the numerical values of $a$, $\nu$ and $T$ are known, it is feasible to evaluate many more coefficients $\al_{2k}$ by this method; see Section 4 for an example.

From (\ref{e25a}) and (\ref{e27a}), we have when $\ta=\pi$
\bee\label{e31a}
f(u)=A\bl\{\frac{1}{u-u_0}+\sum_{r=0}^\infty B_r(\pi) u^r\br\}=A\bl\{-\frac{1}{u_0}\bl(1+\frac{u}{u_0}+\frac{u^2}{u_0^2}+ \ldots\br)+\sum_{r=0}^\infty B_r(\pi) u^r\br\}
\quad (u<u_0),
\ee
whence it follows that the coefficients $B_{2k}(\pi)$ are given by
\bee\label{e31}
B_{2k}(\pi)=ae^{-i\nu\om(\pi)}\al_{2k}+\frac{1}{(ic(\pi))^{2k+1}}.
\ee
The leading coefficient consequently has the value
\bee\label{e32}
B_0(\pi)=\frac{ae^{-i\nu\om(\pi)}}{1-e^{ia\om(\pi)}}+\frac{1}{ic(\pi)}.
\ee

The form (\ref{e31}) has the inconvenient feature of a removable singularity since $\om(\pi)$ and $c(\pi)\to 0$ as $a\to 1$. For
$a\simeq 1$, we can expand the coefficients $B_{2k}(\pi)$ in ascending powers of $\om(\pi)$, viz.
\bee\label{e33}
B_{2k}(\pi)=\sum_{r=0}^\infty b_{2k,r}\,\om^r, \qquad \om\equiv \om(\pi).
\ee
Using {\it Mathematica} to carry out this procedure and the expansion of $c(\pi)$ in powers of $\om(\pi)$ in (\ref{e27}), we obtain the first few values of the coefficients $b_{2k,r}$ in the form:
\begin{eqnarray}
b_{0,0}&=&\frac{1}{2} a+\nu-\frac{1}{6},\quad b_{0,1}=-\frac{i}{12}(a^2+6a\nu+6\nu^2),\nonumber\\ b_{0,2}&=&-\frac{1}{1080}(1+90\nu(a+\nu)(a+2\nu)),\quad
b_{0,3}=\frac{i}{12960}(1+18a^4-540\nu^2(a+\nu)^2), \label{e33a} \\
b_{0,4}&=&\frac{1}{181440}(-1-252a^4\nu+2520a^2\nu^3+3780a\nu^4+1512\nu^5),\ \ldots\, ,\nonumber
\end{eqnarray}
\begin{eqnarray}
b_{2,0}&=&\frac{1}{1080}(-2 + 45 a - 45 a^2 + 90 \nu - 270 a \nu + 90 a^2 \nu - 270 \nu^2 + 
 270 a \nu^2 + 180 \nu^3),\nonumber\\
 b_{2,1}&=&\frac{i}{1440} (-1 - 10 a^2 + 6 a^4 - 60 a \nu + 120 a^2 \nu - 60 \nu^2 + 
   360 a \nu^2 - 180 a^2 \nu^2 + 240 \nu^3\nonumber\\
   &&\hspace{8cm} - 360 a \nu^3 - 180 \nu^4),\label{e33b}\\
b_{2,2}&=&\frac{1}{60480} (1 - 126 a^4 - 420 a^2 \nu + 504 a^4 \nu - 1260 a \nu^2 + 
   3780 a^2 \nu^2 - 840 \nu^3 \nonumber\\
   &&\hspace{3cm}+ 7560 a \nu^3 - 5040 a^2 \nu^3 + 3780 \nu^4 - 
   7560 a \nu^4 - 3024 \nu^5),\ \ldots\  \nonumber
\end{eqnarray}
and
\begin{eqnarray}
b_{4,0}&=&\frac{1}{181440}(65 + 105 a - 630 a^2 + 210 a^4 + 210 \nu - 3780 a \nu + 
  4200 a^2 \nu - 252 a^4 \nu - 3780 \nu^2\nonumber\\
  && + 12600 a \nu^2 - 6300 a^2 \nu^2 + 
  8400 \nu^3 - 12600 a \nu^3 + 2520 a^2 \nu^3 - 6300 \nu^4 + 3780 a \nu^4 + 
  1512 \nu^5),\nonumber\\
b_{4,1}&=&\frac{i}{1088640} (2 + 105 a^2 - 1260 a^4 + 180 a^6 + 630 a \nu - 7560 a^2 \nu + 
   5040 a^4 \nu + 630 \nu^2 \label{e33c}\\
   &&- 22680 a \nu^2 + 37800 a^2 \nu^2 - 
   3780 a^4 \nu^2 - 15120 \nu^3 + 75600 a \nu^3 - 50400 a^2 \nu^3 + 
   37800 \nu^4\nonumber\\
   &&\hspace{3cm} - 75600 a \nu^4 + 18900 a^2 \nu^4 - 30240 \nu^5 + 
   22680 a \nu^5 + 7560 \nu^6),\ \ldots\,.\nonumber
\end{eqnarray}
It may be observed that the even-order coefficients $b_{2k,2r}$ are real whereas the odd-order coefficients $b_{2k,2r+1}$ are imaginary.
\vspace{0.5cm}

\noindent 3.2\ \ {\bf Approximate form of $R_M(a;-x)$ as $a\to 1$}
\vspace{0.3cm}

\noindent
An estimate of the value of the exponentially small term $R_M(a;-x)$ as $a\to1$ can be obtained from Theorem 1 and the fact that, from (\ref{e27}), $c(\pi)=\om+O(\om^2)$. 
%
Then, from (\ref{e210}), it follows that to leading order
\[R_M(a;-x)\simeq \frac{1}{a}\,e^{X\cos \pi/a}\mbox{erfc}\,[\om{\sqrt{X/2}}]-\frac{2e^{-X}}{a\sqrt{2\pi X}}\,
\Re \{ie^{-i\om X} B_0(\pi)\}.\]
Since $B_0(\pi)=b_{0,0}+b_{0,1}\om+O(\om^2)$, where $b_{0,1}$ is imaginary and $\om=\pi(1-a)/a$, we finally obtain 
\[
R_M(a;-x)\simeq\frac{1}{a}\,e^{X\cos \pi/a}\mbox{erfc}\,\bl[\frac{\pi(1-a)}{a} \sqrt{\frac{X}{2}}\br]\hspace{4cm}\]
\bee\label{e310}
\hspace{4cm}-\frac{2e^{-X}}{a\sqrt{2\pi X}}\,\bl\{ b_{0,0} \sin\,\om X+|b_{0,1}|\,\om \cos\,\om X\br\}
\ee
as $x\to+\infty$, where we recall that $X=x^{1/a}$ and $b_{0,0}$, $b_{0,1}$ are given in (\ref{e33a}).

The behaviour of $R_M(a;-x)$ as one approaches the limit $a=1$ is seen to be controlled by a complementary error function, which increases rapidly as $a\to1$ to the value  erfc$(0)=1$ since erfc\,$\xi\sim \exp [-\xi^2]/\sqrt{\pi}\xi$ for $\xi\gg1$. When $a=1$, the quantity $\om\equiv\om(\pi)=0$ and we recover the limiting value of the exponentially small term $R_M(a;-x)=e^{-x}$. Thus the formula (\ref{e310}) correctly describes the appearance of the exponential $e^{-x}$ when $a=1$.

\vspace{0.6cm}

\begin{center}
{\bf 4. \ Numerical results and concluding remarks}
\end{center}
\setcounter{section}{4}
\setcounter{equation}{0}
\renewcommand{\theequation}{\arabic{section}.\arabic{equation}}
To verify the accuracy of the expansion in Theorem 1 we subtract the optimally truncated algebraic expansion
from $E_a(-x)$ and define
\bee\label{e41}
{\cal E}(a;x):=E_a(-x)+\sum_{k=1}^M\frac{(-x)^{-k}}{\g(1-ak)}, 
\ee
where the optimal index $M$ is defined in (\ref{e22}). This quantity is then compared to the exponentially small
contribution $R_M(a;-x)$ for different $a$ and $x$.

\begin{table}[h]
\caption{\footnotesize{The coefficients $B_{2k}(\pi)$ when $a=0.99$ and $x=40$ ($M=42$).}}
\begin{center}
\begin{tabular}{|l|c|}
\hline
&\\[-0.25cm]
\mcol{1}{|c|}{$k$} & \mcol{1}{c|}{$B_{2k}(\pi)$}\\
[.1cm]\hline
&\\[-0.25cm]
0 & $+3.8975364113\times10^{-1} - 3.6166205223\times10^{-3}i$ \\
1 & $-6.4791569264\times10^{-3} - 2.2873163550\times10^{-5}i$ \\
2 & $+1.1939771912\times10^{-3} + 2.9428888000\times10^{-5}i$ \\
3 & $+6.7326294689\times10^{-5} - 3.3561255923\times10^{-7}i$ \\
4 & $+6.4497172230\times10^{-6} - 2.2913466614\times10^{-7}i$ \\
5 & $-4.9612005443\times10^{-7} + 4.0896790580\times10^{-9}i$ \\
6 & $-3.8100530725\times10^{-8} + 1.6905896799\times10^{-9}i$ \\ 
[.1cm]\hline
\end{tabular}
\end{center}
\end{table}

To illustrate we consider the case $x=40$ and $a=0.99$. From (\ref{e22}) we find $M=42$ with the paramneter $\nu=0.0614272718\ldots\ $. The first three coefficients $B_{2k}(\pi)$ can be computed from (\ref{e31}) using the values of $\al_0$, $\al_2$ and $\al_4$ stated in Section 3.1. However, since we have numerical values the higher coefficients can be obtained by the approach discussed at the end of Section 3.1, whereby we expand $f(u)$ in (\ref{e25a}) using the Series command in {\it Mathematica} together with (\ref{e31}). The values of $B_{2k}(\pi)$ for $0\leq k\leq 6$ so obtained are presented in Table 1. We note that the values of $B_{2k}(\pi)$ depend on $x$ through the quantity $\nu$ defined in (\ref{e22}). In Table 2 we show the values of $R_M(a;-x)$ for different truncation index $k$ and two values of $a\simeq1$ compared with the computed values of ${\cal E}(a;x)$. It is seen that there is excellent agreement between the computed value of ${\cal E}(a;x)$ and the asymptotic estimate for the exponentially small contribution.
\begin{table}[h]
\caption{\footnotesize{The values of $R_M(a;-x)$ for different truncation index $k$ when $x=40$: (i) $a=0.99$, $M=42$
 and (ii) $a=0.995$, $M=20$. The final row gives the values of ${\cal E}(a;x)$ defined in (\ref{e41}) for comparison.}}
\begin{center}
\begin{tabular}{|c|c|c|}
\hline
&\\[-0.25cm]
\mcol{1}{|c|}{$k$} & \mcol{1}{c|}{$R_M(a;-x),\ \ a=0.99$} & \mcol{1}{c|}{$R_M(a;-x),\ \ a=0.995$}\\
[.1cm]\hline
&&\\[-0.25cm]
0 & $1.56{\bf 8}95\,52145\,63456\times10^{-19}$ & $1.3789{\bf 9}\,77500\,62528\times10^{-09}$\\
1 & $1.56913\,{\bf 0}8832\,53406\times10^{-19}$ & $1.3789{\bf 1}\,00449\,63445\times10^{-09}$\\
2 & $1.56913\,32{\bf 3}94\,39717\times10^{-19}$ & $1.37890\,9{\bf 8}868\,81488\times10^{-09}$\\
3 & $1.56913\,3223{\bf 5}\,20415\times10^{-19}$ & $1.37890\,9908{\bf 4}\,34786\times10^{-09}$\\
4 & $1.56913\,32232\,6{\bf 1}265\times10^{-19}$ & $1.37890\,99085\,{\bf 2}9609\times10^{-09}$\\
5 & $1.56913\,32232\,65{\bf 5}55\times10^{-19}$ & $1.37890\,99085\,08{\bf 3}09\times10^{-09}$\\
6 & $1.56913\,32232\,6564{\bf 4}\times10^{-19}$ & $1.37890\,99085\,081{\bf 4}4\times10^{-09}$\\
[.1cm]\hline
&&\\[-0.25cm]
${\cal E}(a;x)$ & $1.56913\,32232\,65642\times10^{-19}$ & $1.37890\,99085\,08192\times10^{-09}$\\
[.1cm]\hline
\end{tabular}
\end{center}
\end{table}

In Table 3 we show values of ${\cal E}(a;x)$ and $R_M(a;-x)$ (with truncation index $k=5$) for a range of $a$-values.
At the end of Section 2.1 it was argued that the parameter $a>\f{1}{3}$ for the sector of validity of the expansion (\ref{e210}) to include the negative real axis. It is noteworthy that there continues to be good agreement between ${\cal E}(a;x)$ and $R_M(a;-x)$ even when $a\leq\f{1}{3}$. The value of $x$ chosen in the cases $a=\f{1}{3}$ and $a=\f{1}{4}$ in Table 3 is small; larger values would result in very large optimal truncation index $M$ (for example, if $x=10$ when $a=\f{1}{4}$, we find $M=4\times 10^4$). This would produce extreme accuracy from just the algebraic expansion, with $R_M(a;-x)$ so small as to be negligible in most applications.
The validity of this agreement when $a\leq\f{1}{3}$ would require further investigation, which is not carried out here as our main interest is in the limit $a\to1$.
\begin{table}[th]
\caption{\footnotesize{The values of ${\cal E}(a;x)$ and $R_M(a;-x)$ (with truncation index $k=5$) for different values of $a$ and $x$.}}
\begin{center}
\begin{tabular}{|crr|c|c|}
\hline
&&&&\\[-0.25cm]
\mcol{1}{|c}{$a$} & \mcol{1}{c}{$x$} & \mcol{1}{c}{$M$} & \mcol{1}{|c|}{${\cal E}(a;x)$} & \mcol{1}{c|}{$R_M(a;-x)$}\\
[.1cm]\hline
&&&&\\[-0.25cm]
0.95 & 20 & 25 & $-2.521343\,284521\times 10^{-11}$ & $-2.521343\,28452{\bf 2}\times 10^{-11}$\\
0.90 & 20 & 21 & $-2.706560\,459479\times 10^{-13}$ & $-2.706560\,45947{\bf 8}\times 10^{-13}$\\
0.80 & 20 & 53 & $-4.827618\,810882\times 10^{-20}$ & $-4.827618\,810882\times 10^{-20}$\\
0.70 & 15 & 68 & $-3.052228\,407002\times 10^{-23}$ & $-3.052228\,407002\times 10^{-23}$\\
0.60 & 10 & 77 & $-6.895973\,422484\times 10^{-22}$ & $-6.895973\,422484\times 10^{-22}$\\
0.50 & 5  & 50 & $-1.106145\,146730\times 10^{-12}$ & $-1.106145\,146730\times 10^{-12}$\\
0.3${\hat 3}$ & 3 & 81 & $+8.345377\,837784\times 10^{-14}$ & $+8.345377\,8377{\bf 3}5\times 10^{-14}$\\
0.25 & 3  & 324 & $-1.220075\,244872\times 10^{-37}$ & $-1.220075\,244872\times10^{-37}$\\
[.1cm]\hline
\end{tabular}
\end{center}
\end{table}

\vspace{0.6cm}

\begin{center}
{\bf Appendix: Estimation of the contribution along the upper side of the cut}
\end{center}
\setcounter{section}{1}
\setcounter{equation}{0}
\renewcommand{\theequation}{\Alph{section}.\arabic{equation}}
The procedure we employ is a slight modification of that described by Olver \cite{Olv}
in the treatment of the generalised exponential integral; see also \cite[\S 6.2.6]{PK}.
If we make the change of variable $t=e^{-\pi i}\tau$ in (\ref{e22a}), the integral taken along the upper side of the cut in the $\tau$-plane becomes
\bee\label{e23}
J=e^{-iM(\ta-\pi a)}\frac{e^{-X}}{2\pi i}\int_0^\infty e^{-X\psi(t)}\,\frac{t^{a+\nu-1}}{t^a-t_0^a}\,dt,
\ee
where
\[\psi(t)=t-\log\,t-1,\quad t_0=e^{i\om(\ta)},\quad \om(\ta)=(\ta-\pi a)/a.\]
In the $t$-plane, the branch cut is now situated on $[0,\infty)$ and
the integration path in (\ref{e23}) passes {\it below} the image of the pole $P_0$ and the saddle at $t=1$.
Setting
\bee\label{e24}
\fs u^2=t-\log\,t-1,\qquad \frac{dt}{du}=\frac{ut}{t-1},
\ee
we can express the integral (\ref{e23}) in the form
\bee\label{e25}
J=e^{-iM(\ta-\pi a)}\, \frac{e^{-X}}{2\pi i}\int_{-\infty}^\infty e^{-\fr xu^2} f(u)\,du,
\ee
where
\bee\label{e25b}
f(u):=\frac{t^{a+\nu-1}}{t^a-t_0^a}\,\frac{dt}{du}=\frac{ut^{a+\nu}}{(t-1)(t^a-t_0^a)}.
\ee

The function $f(u)$ can be expanded in the form
\bee\label{e25a}
f(u)=A\bl\{\frac{1}{u-u_0}+g(u)\br\},
\ee
where the pole at $u=u_0\equiv ic(\ta)$ corresponds to the pole in the $t$-plane at $t_0=e^{i\om(\ta)}$ and $g(u)$ is analytic at the point $u=u_0$. We have from (\ref{e24})
\bee\label{e26}
\fs c^2(\ta)=1+i\om(\ta)-e^{i\om(\ta)},
\ee
where the branch of the square root is chosen so that near $\ta=\pi a$ the expansion of $c(\ta)$ has the form
\bee\label{e27}
c(\ta)=\om(\ta)+\f{1}{6}i\om^2(\ta)-\f{1}{36}\om^3(\ta)-\f{1}{270}i\om^4(\ta)+\f{1}{2592}\om^5(\ta)+\ldots\ .
\ee
The constant $A$ appearing in (\ref{e25a}) can be determined by a limiting process. If we let $t=t_0+\epsilon$, $\epsilon\to0$, so that from (\ref{e24}) $u-u_0=\epsilon(t_0-1)/(u_0 t_0)+O(\epsilon^2)$, we find 
\[A=\mathop{\lim_{\scriptstyle u\to u_0 \atop\scriptstyle t\to t_0}} (u-u_0)f(u)=\frac{e^{i\nu\om(\ta)}}{a}.\]

Substitution of the above expansion for $f(u)$ in (\ref{e25}) then yields
\[J=\frac{e^{-X-i\om(\ta)X}}{2\pi i}\bl\{\int_{-\infty}^\infty \frac{e^{-\fr Xu^2}}{u-u_0}\,du+\int_{-\infty}^\infty e^{-\fr Xu^2}g(u)\,du\br\}.\]
The first integral on the right-hand side of the above expression (where the path is indented to pass below the pole $u_0$) can be evaluated in terms of the complementary error function
\[\int_{-\infty}^\infty \frac{e^{-\fr Xu^2}}{u-u_0}\,du=\pi i \,e^{\fr Xc^2(\ta)} \mbox{erfc}\,[c(\ta) \sqrt{X/2}\,].\]
In the second integral the path may be taken as the real axis with no indentation, since the integrand has no singularity on the integration path. If we expand $g(u)$ as a Maclaurin series
\bee\label{e27a}
g(u)=\sum_{r=0}^\infty B_r(\ta) u^r,
\ee
we find
\[\int_{-\infty}^\infty e^{-\fr Xu^2}g(u)\,du \sim \pi^{1/2} \sum_{k=0}^\infty B_{2k}(\ta) (\fs)_k (\fs X)^{-k}\qquad (X\to\infty).\]

Collecting together these results and noting that $$e^{\fr Xc^2(\ta)}=\exp\,[z^{1/a}+X+i\om(\ta) X]$$ by (\ref{e26}), we finally obtain from Theroem 1 of \cite[p.~1473]{Olv} the desired expansion\footnote{There is an error in the sign of the second term in this expansion in \cite[(2.8)]{P}.} 
\bee\label{e28}
J\sim \frac{1}{a}\bl\{\frac{\exp\,[z^{1/a}]}{2}\,\mbox{erfc}\,[c(\ta) \sqrt{X/2}\,]-\frac{ie^{-X-i\om(\ta) X}}{\sqrt{2\pi X}}\sum_{k=0}^\infty B_{2k}(\ta) (\fs)_k (\fs X)^{-k}\br\}
\ee
as $|z|\to\infty$ in the sector $-\pi a<\theta<3\pi a$,
where $c(\ta)$ is defined by (\ref{e26}) with the expansion in ascending powers of $\om(\ta)$ given in (\ref{e27}). The coefficients $B_{2k}(\ta)\equiv B_{2k}(\ta,\nu)$ in the case $\ta=\pi$ are discussed in Section 3.2. The above sector clearly includes the negative real axis $\arg\,z=\pi$ when $a>\f{1}{3}$.

\end{document}